\documentstyle[12pt, fullpage]{amsart}
\newtheorem{lemma}{Lemma}

\newenvironment{pflike}[1]{\noindent{\bf #1}}{\vskip10pt} 
\newenvironment{proof}{\begin{pflike}{Proof:}}{\qed\end{pflike}}

\begin{document}
\thispagestyle{empty}

\title{Absolutely Abnormal Numbers}
\author{Greg Martin}
\address{Department of Mathematics\\University of Toronto\\Canada M5S 3G3}
\email{gerg@@math.toronto.edu}
\maketitle

\section{Introduction}

A normal number is one whose decimal expansion (or expansion
to some base other than~10) contains all possible finite
configurations of digits with roughly their expected frequencies. More
formally, when $b\ge2$ is an integer, let
\begin{equation}
N(\alpha;b,a,x) = \#\{ 1\le n\le x\colon \text{the $n$th digit
in the base-$b$ expansion of $\alpha$ is }a\}
\label{Ndef}
\end{equation}
denote the counting function of the occurrences of the digit $a$
($0\le a<b$) in the $b$-ary expansion of the real number $\alpha$, and
define the corresponding limiting frequency
\begin{equation}
\textstyle \delta(\alpha;b,a) = \lim\limits_{x\to\infty}
x^{-1}N(\alpha;b,a,x),
\label{deltadef}
\end{equation}
if the limit exists. The number $\alpha$ is {\it simply normal\/} to
the base $b$ if the limit defining $\delta(\alpha;b,a)$ exists and
equals $1/b$ for each $0\le a<b$. (When $\alpha$ is a $b$-adic
fraction $a/b^n$, which has one $b$-ary expansion with all but
finitely many digits equaling zero and another $b$-ary expansion with
all but finitely many digits equaling $b-1$, these limiting
frequencies are not uniquely defined; however, such an $\alpha$ will
not be simply normal to the base $b$ in either case.) A number is {\it
normal\/} to the base $b$ if it is simply normal to each of the bases
$b$, $b^2$, $b^3$, \dots. This is equivalent (see \cite[Chapter
8]{niven}) to demanding that for any finite string $b_1b_2\dots b_k$
of base-$b$ digits, the limiting frequency of occurrences of this
string in the $b$-ary expansion of $\alpha$ (defined analogously to
equation~(\ref{deltadef}) above) exists and equals $1/b^k$.

For instance, it was shown by Champernowne~\cite{champ} that the
number $0.12345678910111213\dots$ formed by concatenating all of the
positive integers together into a single decimal is normal to base 10
(the analogous construction works for any base $b\ge2$), and this sort
of example has been generalized~\cite{CE,DE}. It is known that almost
all real numbers are normal to any given base $b$ (see for
instance~\cite[Theorem 8.11]{niven}), and consequently almost all real
numbers are {\it absolutely normal\/}, i.e., normal to all bases
$b\ge2$ simultaneously. On the other hand, we have not proven a single
naturally occurring real number to be absolutely normal.

Let us call a number {\it abnormal\/} to the base $b$ if it is not
normal to the base $b$, and {\it absolutely abnormal\/} if it is
abnormal to all bases $b\ge2$ simultaneously. For instance, every
rational number $r$ is absolutely abnormal: any $b$-ary expansion of
$r$ will eventually repeat, say with period $k$, in which case $r$ is
about as far from being simply normal to the base $b^k$ as it can
be. Even though the set of absolutely abnormal numbers is the
intersection of countably many sets of measure zero, it was pointed
out by Maxfield~\cite{max} that the set of numbers normal to a given
base $b$ is uncountable and dense; later, Schmidt~\cite{schmidt} gave
a complicated constructive proof of this fact. In this paper we
exhibit a simple construction of a specific irrational (in fact,
transcendental) real number that is absolutely abnormal:

\medskip\noindent{\bf Theorem.} {\it The number $\alpha$ defined
in equation~(\ref{alphaprod}) below is irrational and absolutely
abnormal.}

\noindent In fact, our construction easily generalizes to
produce concretely an uncountable set of absolutely abnormal numbers
in any open interval.

It is instructive to consider why constructing an irrational,
absolutely abnormal number is even difficult. Since we already know
that rational numbers are absolutely abnormal, our first thought might
be to choose an irrational number whose $b$-ary expansions mimic those
of rational numbers for long stretches, i.e., an irrational number
with very good rational approximations. Thus a natural class to
consider is the {\it Liouville numbers\/}, defined to be those real
numbers $\beta$ such that for every positive integer $m$, there exists
a rational number $\frac pq$ (not necessarily in lowest terms) satisfying
\begin{equation}
0 < \Big| \beta-\frac pq \Big| < {1\over q^m}.  \label{Liouproperty}
\end{equation}
These Liouville numbers are all transcendcental (see
Lemma~\ref{Lioulem} below)---in fact Liouville introduced these
numbers precisely to exhibit specific transcendental numbers, and the
off-cited example
\begin{equation*}
\beta = \sum_{n=1}^\infty 10^{-n!} = 0.11000100000000000000000100\dots
\end{equation*}
is usually the first number that students see proven transcendental.

Clearly $\beta$ is abnormal to the base 10; how would we go about
showing, for example, that $\beta$ is abnormal to the base 2? We would
try to argue that the binary expansion of $\beta$ agrees with that of
each of the rational numbers
\begin{equation}
\beta_k = \sum_{n=1}^k 10^{-n!}  \label{Lioufinite}
\end{equation}
through about the $((n+1)!\log_210)$-th binary digit. Since
each $\beta_k$ is rational and thus abnormal to the base 2, can we
conclude that $\beta$ itself is abnormal to the base 2?

Not quite: it seems that we would have to show that there is a fixed
power $2^n$ such that infinitely many of the $\beta_k$ were not simply
normal to the base $2^n$. (For each $\beta_k$ there is {\it some\/} power
$2^{n_k}$ such that $\beta_k$ is not simply normal to the base
$2^{n_k}$, but these exponents $n_k$ might very well grow with $k$.)
In fact, it is not hard to show (using the fact that 2 is a primitive
root modulo every power of 5) that any 10-adic fraction that is not a
2-adic fraction---including each $\beta_k$---{\it is\/} simply normal to
the base 2! In general, without actually computing binary expansions of
specific fractions, it seems impossible to rule out the incredible
possibility that the $\beta_k$ are accidentally simply normal to
bases that are high powers of 2. In summary, any Liouville number we write
down is almost certain (morally) to be absolutely abnormal, but actually
proving its absolute abnormality is another matter.

To circumvent this difficulty, we construct a Liouville number whose
successive rational approximations are $b$-adic fractions with $b$
varying, rather than all being 10-adic fractions as in
equation~(\ref{Lioufinite}). The existence of such Liouville numbers
can certainly be proven using just the fact that the $b$-adic
fractions are dense for any integer $b\ge2$; however, our construction
is completely explicit. We first give the complete construction of our
irrational, absolutely abnormal number $\alpha$ and then show
afterwards that $\alpha$ has the required properties.

\section{The construction and proof}

We begin by defining a sequence of integers
\begin{equation*}
d_2=2^2, \quad d_3=3^2, \quad d_4=4^3, \quad d_5 = 5^{16}, \quad d_6
= 6^{\text{30,517,578,125}}, \,\dots
\end{equation*}
with the recursive rule
\begin{equation}
d_j = j^{d_{j-1}/(j-1)} \quad(j\ge3).  \label{djdef}
\end{equation}
This sequence exhibits the pattern
\newlength{\cthree}
\newlength{\ctwo}
\begin{equation*}
d_4=4^{3^{2-1}}, \quad
d_5=5^{4^{\left(
    \setbox3\hbox{$\scriptstyle3^{2-1}-1$}
   \setlength{\cthree}{-0.5\ht3}\addtolength{\cthree}{2pt}
   \raisebox{\cthree}{\usebox3}
   \right)}}, \quad
d_6=6^{5^{\left(\!
   \setbox2\hbox{$\scriptstyle4^{ \left(
    \setbox3\hbox{$\scriptstyle3^{2-1}-1$}
   \setlength{\cthree}{-0.5\ht3}\addtolength{\cthree}{2pt}
   \raisebox{\cthree}{\usebox3}
   \right) } -1 $}
  \setlength{\ctwo}{-0.5\ht2}\addtolength{\ctwo}{2pt}
  \raisebox{\ctwo}{\usebox2}
  \right) }},\, \dots
\end{equation*}
which in general gives the typesetting nightmare
\newlength{\cone}
\newlength{\czero}
\newlength{\cfour}
\newcommand{\adots}{\mathinner{\mkern2mu\raise0.4pt\hbox{.}\mkern2mu\raise1.7pt\hbox{.}\mkern2mu\raise3pt\hbox{.}\mkern1mu}}
\newcommand{\bdots}{\mathinner{\mkern2mu\raise3pt\hbox{.}\mkern2mu\raise1.7pt\hbox{.}\mkern2mu\raise0.4pt\hbox{.}\mkern1mu}}
\begin{equation}
d_j=j^{ \,(j-1)^{ \!\!\!\left(
 \setbox0\hbox{$\scriptstyle(j-2)^{ \!\!\left(
 \setbox4\hbox{$\scriptstyle(j-3)^{ \!\!\left(
  \setbox1\hbox{$\scriptstyle{\adots\vphantom{\big(}}^{ \left(
   \setbox2\hbox{$\scriptstyle4^{ \left(
    \setbox3\hbox{$\scriptstyle3^{2-1}$}
   \setlength{\cthree}{-0.5\ht3}\addtolength{\cthree}{2pt}
   \raisebox{\cthree}{\usebox3}
   \right) } {}-1 $}
  \setlength{\ctwo}{-0.5\ht2}\addtolength{\ctwo}{2pt}
  \raisebox{\ctwo}{\usebox2}
  \right) } \bdots $} 
 \setlength{\cone}{-0.5\ht1}\addtolength{\cone}{4pt}
 \raisebox{\cone}{\usebox1}
 \right) } \!\!{}-1 $}
 \setlength{\cfour}{-0.5\ht4}\addtolength{\cfour}{4pt}
 \raisebox{\cfour}{\usebox4}
 \right) } \!\!{}-1 $}
\setlength{\czero}{-0.5\ht0}\addtolength{\czero}{4pt}
\raisebox{\czero}{\usebox0}
\right) }}.
\label{nightmare}
\end{equation}
Using these integers, we define the sequence of rational numbers
\begin{equation}
\alpha_k = \prod_{j=2}^k \big( 1-\frac1{d_j} \big),  \label{alphakdef}
\end{equation}
so that $\alpha_2=\frac14$, $\alpha_3=\frac23$,
$\alpha_4={21\over32}$, $\alpha_5={\text{100,135,803,222} \over
\text{152,587,890,625}}$, and so on.

We now nominate
\begin{equation}
\alpha = \lim_{k\to\infty} \alpha_k = \prod_{j=2}^\infty
\big(1-\frac1{d_j}\big)
\label{alphaprod}
\end{equation}
as our candidate for an irrational, absolutely abnormal number. The
first few digits in the decimal expansion of $\alpha$ are
\begin{equation}
\alpha = 0.6562499999956991\!
\underbrace{99999\dots99999}_{\text{23,747,291,559 9s}}\!
8528404201690728\dots, 
\label{alphadecimal}
\end{equation}
from which we can get an inkling of the extreme abnormality of
$\alpha$ (at least to the base~10). We need to prove three things
concerning this number $\alpha$: first, that the infinite
product~(\ref{alphaprod}) defining $\alpha$ actually converges; second,
that $\alpha$ is irrational; and finally, that $\alpha$ is absolutely
abnormal.

It is apparent from the expressions~(\ref{djdef})
and~(\ref{nightmare}) that the $d_j$ grow (ridiculously) rapidly and
hence that the infinite product~(\ref{alphaprod}) should indeed
converge. The following lemma provides a crude inequality relating the
integers $d_j$ that we can use to prove this assertion rigorously.

\begin{lemma}
For $j\ge5$ we have $d_j > 2d_{j-1}^2$.  \label{cubelem}
\end{lemma}

\begin{proof}
We proceed by induction, the case $j=5$ being true by inspection. For
$j>5$ we surely have
\begin{equation*}
d_j = j^{d_{j-1}/(j-1)} > 5^{d_{j-1}/(j-1)}.
\end{equation*}
Notice that from the definition~(\ref{djdef}) of the $d_j$,
\begin{equation*}
{d_{j-1}\over j-1} = (j-1)^{{\scriptstyle
d_{j-2}\over\scriptstyle\mathstrut j-2}-1} > (j-1)^{{\scriptstyle
d_{j-2}\over\scriptstyle\mathstrut2(j-2)}} = \sqrt{d_{j-1}},
\end{equation*}
and therefore
\begin{equation*}
d_j > 5^{\sqrt{d_{j-1}}}.
\end{equation*}
Now using the fact (easily proven by your favorite calculus student)
that $5^x \ge x^5$ for $x\ge5$, we conclude that
\begin{equation*}
d_j > \Big(\sqrt{d_{j-1}}\Big)^5 > 2d_{j-1}^2,
\end{equation*}
as desired.
\end{proof}

Equipped with this inequality, we can now show that the infinite
product~(\ref{alphaprod}) defining $\alpha$ converges. Moreover, we
can show that the number $\alpha$ is well approximated by
the rational numbers $\alpha_k$. (Notice that $\alpha_4$ is exactly
$0.65625$ and $\alpha_5$ is exactly $0.6562499999956992$---cf.~the
decimal expansion~(\ref{alphadecimal}) of $\alpha$.)

\begin{lemma}
The product\/~{\rm (\ref{alphaprod})} defining $\alpha$
converges. Moreover, for $k\ge2$ we have
\begin{equation}
\alpha_k > \alpha > \alpha_k - \frac2{d_{k+1}}.  \label{close}
\end{equation}
\label{squeezelem}
\end{lemma}

\begin{proof}
To show that the product~(\ref{alphaprod}) converges, we must show
that the corresponding sum $\sum_{j=2}^\infty 1/d_j$ converges. But by
Lemma \ref{cubelem} we certainly have $d_j>2d_{j-1}$ for $j\ge5$,
and therefore
\begin{equation*}
\sum_{j=2}^\infty \frac1{d_j} \le \frac1{d_2} + \frac1{d_3} + \sum_{j=4}^\infty
\frac1{2^{j-4}d_4} = \frac1{d_2} + \frac1{d_3} + \frac2{d_4} < \infty.
\end{equation*}
Similarly, using the fact that $1\ge\prod(1-x_j) \ge 1-\sum x_j$ for any real
numbers $0\le x_j\le1$, we see that for $k\ge3$
\begin{multline*}
\alpha_k > \alpha = \alpha_k \prod_{j=k+1}^\infty \big( 1-\frac1{d_j} \big) \ge
\alpha_k \bigg( 1-\sum_{j=k+1}^\infty \frac1{d_j} \bigg) \\
> \alpha_k \bigg( 1-\sum_{j=k+1}^\infty \frac1{2^{j-k-1}d_{k+1}} \bigg) =
\alpha_k \big( 1 - \frac2{d_{k+1}} \big) > \alpha_k - \frac2{d_{k+1}}.
\end{multline*}
The inequalities~(\ref{close}) for $k=2$ follow from those for $k=3$,
as it is easily verified by hand that $\alpha_2 > \alpha_3 > \alpha_3
- \frac2{d_4} > \alpha_2 - \frac2{d_3}$.
\end{proof}

It turns out that both the proof that $\alpha$ is irrational and the proof
that $\alpha$ is absolutely abnormal hinge on the fact that each
rational approximation $\alpha_k$ is in fact a $k$-adic fraction---that
is, when $\alpha_k$ is expressed in lowest terms, its denominator divides
a power of $k$. In other words, each time we multiply $\alpha_{k-1}$ by
$1-{1\over d_k} = {d_k-1\over d_k}$ to obtain $\alpha_k$, the numerator of the latter
fraction completely cancels out the denominator of $\alpha_{k-1}$, so
that all that can remain in the denominator of $\alpha_k$ are the powers
of $k$ present in $d_k$. Proving that this always happens is an exercise
in elementary number theory, which we present in the next three lemmas.

\begin{lemma}
Let $k$ and $r$ be positive integers, and let $p$ be a prime. If $k$
is divisible by $p^r$, then $(k+1)^p-1$ is divisible by $p^{r+1}$.
\label{psteplem}
\end{lemma}

\begin{proof}
Writing $k=p^rn$, we have from the binomial theorem
\begin{equation*}
\begin{split}
(k+1)^p-1 &= (p^rn+1)^p-1 \\
&= \big\{ (p^rn)^p + {\textstyle{p\choose p-1}}(p^rn)^{p-1} + \dots +
{\textstyle{p\choose2}}(p^rn)^2 + {\textstyle{p\choose1}}p^rn + 1 \big\} - 1
\\
&= p^{rp}n^p + {\textstyle{p\choose p-1}}p^{r(p-1)}n^{p-1} + \dots +
{\textstyle{p\choose2}}p^{2r}n^2 + p\cdot p^rn.
\end{split}
\end{equation*}
Since all of these binomial coefficients $p\choose k$ are integers,
each term in this last sum is visibly divisible by $p^{r+1}$.
\end{proof}

\begin{lemma}
For any positive integers $k$ and $m$, the integer $(k+1)^{k^m}-1$ is
divisible by $k^{m+1}$.\label{refereelem}
\end{lemma}

\begin{proof}
If $p^r$ is any prime power dividing $k$, an $rm$-fold application of
Lemma~\ref{psteplem} shows us that $(k+1)^{p^{rm}}-1$ is divisible by
$p^{r+rm}$. Then, since
\begin{equation*}
(k+1)^{k^m}-1 = \big( (k+1)^{p^{rm}}-1 \big) \big( (k+1)^{k^m-p^{rm}}
+ (k+1)^{k^m-2p^{rm}} + \dots + (k+1)^{p^{rm}} + 1 \big),
\end{equation*}
we see that $(k+1)^{k^m}-1$ is also divisible by $p^{r(m+1)}$.

In particular, since $p^r$ was an arbitrary prime power dividing $k$,
we see that $(k+1)^{k^m}-1$ is divisible by every prime power that
divides $k^{m+1}$. This is enough to verify that $(k+1)^{k^m}-1$ is
divisible by $k^{m+1}$ itself.
\end{proof}

\begin{lemma}
For each $k\ge2$, the product $d_k\alpha_k$ is an integer.  In
particular, since $d_k$ is a power of $k$, we see that $\alpha_k$ is a
$k$-adic fraction.
\label{intlem}
\end{lemma}

\begin{proof}
We proceed by induction on $k$, the cases $k=2$ and $k=3$ being
evident by inspection. For the inductive step, suppose (as our
induction hypothesis) that $d_k\alpha_k$ is indeed an integer for a
given $k\ge3$. We may write
\begin{equation}
d_{k+1}\alpha_{k+1} = (d_{k+1}-1)\alpha_k = \big( (k+1)^{d_k/k}-1
\big) \alpha_k = {(k+1)^{d_k/k}-1\over d_k} \cdot d_k\alpha_k
\label{integer}
\end{equation}
by the definitions~(\ref{djdef}) and~(\ref{alphakdef}) of $d_k$ and
$\alpha_k$, respectively. The second factor $d_k\alpha_k$ is an
integer by the induction hypothesis. On the other hand, we may rewrite
\begin{equation*}
(k+1)^{d_k/k}-1 = (k+1)^{k^{d_{k-1}/(k-1)-1}}-1.
\end{equation*}
Applying Lemma~\ref{refereelem} with $m=d_{k-1}/(k-1)-1$, we see that
this expression is divisible by $k^{d_{k-1}/(k-1)} = d_k$. Therefore
the fraction on the right-hand side of equation~(\ref{integer}) is
in fact an integer, and so $d_{k+1}\alpha_{k+1}$ is itself an
integer, which completes the proof.
\end{proof}

As mentioned in the introduction, the key to proving that $\alpha$ is
irrational is to show that it is in fact a Liouville number. It is a
standard fact that any Liouville number is transcendental (see for
instance~\cite[Theorem 7.9]{niven}); for the sake of keeping this
paper self-contained, we include a proof.

\begin{lemma}
Every Liouville number is transcendental.  \label{Lioulem}
\end{lemma}

\begin{proof}
We shall prove the contrapositive, that no algebraic number can satisfy
the Liouville property~(\ref{Liouproperty}) for all positive $m$. Suppose
that $\beta$ is algebraic. Without loss of generality, we may suppose that
$|\beta|\le\frac12$ by adding an appropriate integer. Let
\begin{equation*}
m_\beta(x) = {c_d}x^d + c_{d-1}x^{d-1} + \dots + c_2x^2 + c_1x + c_0
\end{equation*}
be the minimal polynomial for $\beta$, where the coefficients $c_i$
are integers. Now suppose that $\frac pq$ is a rational approximation
to $\beta$, say $|\beta-\frac pq|<\frac12$. Then
\begin{equation*}
\begin{split}
m_\beta\big(\frac pq\big) &= m_\beta\big(\frac pq\big) -
m_\beta(\beta) \\
&= c_d\Big( \big(\frac pq\big)^d - \beta^d \Big) + \dots + c_2\Big(
\big(\frac pq\big)^2 - \beta^2 \Big) + c_1\big(\frac pq-\beta\big) \\
&= \big(\frac pq-\beta\big) \bigg( c_d \Big( \big(\frac pq\big)^{d-1}
+ \big(\frac pq\big)^{d-2}\beta + \dots + \beta^{d-1} \Big) + \dots +
c_2 \Big( \frac pq + \beta \Big) + c_1 \bigg).
\end{split}
\end{equation*}
Since neither $\beta$ nor $\frac pq$ exceeds 1 in absolute value, we
see that
\begin{equation}
\Big| m_\beta\big(\frac pq\big) \Big| \le \Big| \frac pq-\beta \Big|
C(\beta), \label{lowerbeta}
\end{equation}
where we have defined the constant
\begin{equation*}
C(\beta) = d|c_d| + (d-1)|c_{d-1}| + \dots + 2|c_2| + |c_1|.
\end{equation*}

On the other hand, $m_\beta\big(\frac pq\big)$ is a rational number
with denominator at most $q^d$, and it is nonzero since $m_\beta$ is
irreducible. Therefore
\begin{equation}
\Big| m_\beta\big(\frac pq\big) \Big| \ge {1\over q^d}.  \label{upperbeta}
\end{equation}
Together, the inequalities~(\ref{lowerbeta}) and~(\ref{upperbeta})
imply that
\begin{equation*}
\Big| \beta-\frac pq \Big| \ge {C(\beta)^{-1}\over q^d},
\end{equation*}
which precludes the inequality~(\ref{Liouproperty}) from holding when $m$
is large enough.
\end{proof}

To show that $\alpha$ is indeed a Liouville number, we will need an
inequality somewhat stronger than the one given in Lemma~\ref{cubelem}.
The following lemma furnishes a simple inequality that is strong enough
for this purpose.

\begin{lemma}
For $j\ge5$ we have $d_{j+1} > d_j^{\,d_{j-1}}$.  \label{explem}
\end{lemma}

\begin{proof}
It is immediate that
\begin{equation*}
d_{j+1} = (j+1)^{d_j/j} > j^{d_j/j} > j^{2d_{j-1}^2/j} =
(j^{d_{j-1}/(j-1)})^{d_{j-1}\cdot2(j-1)/j} =
d_j^{\,d_{j-1}\cdot2(j-1)/j} > d_j^{\,d_{j-1}},
\end{equation*}
where we have used Lemma~\ref{cubelem} for the second inequality.
\end{proof}

\begin{lemma}
$\alpha$ is a Liouville number; in particular, $\alpha$ is transcendental.
\label{alphaLioulem}
\end{lemma}

\begin{proof}
We can easily show show that the $\alpha_k$ provide the very close
rational approximations needed in equation~(\ref{Liouproperty}) to make
$\alpha$ a Liouville number. Indeed, $\alpha_k$ can be written as a
fraction whose denominator is $d_k$ by Lemma \ref{intlem}, while
Lemma~\ref{squeezelem} tells us that for $k\ge5$
\begin{equation*}
0< |\alpha-\alpha_k| < \frac2{d_{k+1}} < \frac2{d_k^{\,d_{k-1}}},
\end{equation*}
where the last inequality is by Lemma~\ref{explem}. Since $d_{k-1}$
tends to infinity with $k$, this shows that $\alpha$ is a
Liouville number (and hence transcendental by Lemma~\ref{Lioulem}).
\end{proof}

\goodbreak At last we have all the tools we need to establish the theorem
stated in the introduction:

\medskip\noindent{\bf Theorem.} {\it The number $\alpha$ defined
in equation~(\ref{alphaprod}) is irrational and absolutely abnormal.}
\medskip

\begin{proof}
We have just shown in Lemma~\ref{alphaLioulem} that $\alpha$ is
irrational. As for proving that $\alpha$ is absolutely abnormal, the idea
is that for every integer base $b\ge2$, the number $\alpha$ is just a tiny
bit less than the $b$-adic fraction $\alpha_b$. Since the $b$-ary
expansion of $\alpha_b$ terminates in an infinite string of zeros, the
slightly smaller number $\alpha$ will have a long string of digits equal
to $b-1$ before resuming a more random behavior. (This is evident in the
decimal expansion~(\ref{alphadecimal}) of $\alpha$, as $\alpha_5$ is a
10-adic fraction as well as a 5-adic fraction.) This happens more than
once, as each of $\alpha_b$, $\alpha_{b^2}$, $\alpha_{b^3}$, and so on is
a $b$-adic fraction. Consequently, the $b$-ary expansion of $\alpha$ will
have increasingly long strings consisting solely of the digit $b-1$,
which will prevent it from being even simply normal to the base $b$.

More quantitatively, let $b\ge2$ and $r$ be positive integers. Since
$d_{b^r}\alpha_{b^r}$ is an integer by Lemma~\ref{intlem}, and since
$d_{b^r}=(b^r)^{d_{b^r-1}/(b^r-1)}$ by definition, the $b$-ary
expansion of $\alpha_{b^r}$ terminates after at most
$rd_{b^r-1}/(b^r-1)$ nonzero digits. On the other hand, by
Lemma~\ref{squeezelem} we know that $\alpha$ is less than
$\alpha_{b^r}$ but by no more than $2/d_{b^r+1} =
2/(b^r+1)^{d_{b^r}/b^r} < 2/b^{rd_{b^r}/b^r}$. Therefore, when we
subtract this small difference from $\alpha_{b^r}$, the resulting
$b$-ary expansion will have occurrences of the digit $b-1$ beginning
at the $(rd_{b^r-1}/(b^r-1)+1)$-th digit at the latest, and continuing
through at least the $(rd_{b^r}/b^r-1)$-th digit since the difference
will start to show only in the $(rd_{b^r}/b^r)$-th digit at the
soonest.  Using the notation defined in equation~(\ref{Ndef}), this
implies that
\begin{equation*}
N\big(\alpha;b,b-1,{rd_{b^r}\over b^r}\big) \ge {rd_{b^r}\over b^r}
-{rd_{b^r-1}\over b^r-1} -1 > {rd_{b^r}\over b^r} - {2rd_{b^r-1}\over b^r}.
\end{equation*}

At this point we can calculate that
\begin{equation*}
\begin{split}
\limsup_{x\to\infty} x^{-1} N(\alpha;b,b-1,x) &\ge
\limsup_{r\to\infty} \big( {b^r\over rd_{b^r}}
N\big(\alpha;b,b-1,{rd_{b^r}\over b^r}\big) \big) \\
&\ge \limsup_{r\to\infty} \big( 1 - {2d_{b^r-1}\over d_{b^r}} \big).
\end{split}
\end{equation*}
Using Lemma~\ref{cubelem}, we see that
\begin{equation*}
\limsup_{x\to\infty} x^{-1} N(\alpha;b,b-1,x) \ge \limsup_{r\to\infty}
\big( 1 - {2d_{b^r-1}\over 2d_{b^r-1}^2} \big) = \limsup_{r\to\infty}
\big( 1 - {1\over d_{b^r-1}} \big) = 1.
\end{equation*}
In particular, the frequency $\delta(\alpha;b,b-1)$ defined in
equation~(\ref{deltadef}) either does not exist or else equals 1, either of
which precludes $\alpha$ from being simply normal to the base $b$. Since
$b\ge2$ was arbitrary, this shows that
$\alpha$ is absolutely abnormal.
\end{proof}

\section{Generalizations and further questions}

We mentioned in the introduction that our construction of an
irrational, absolutely abnormal number can be generalized to exhibit
an uncountable set of absolutely abnormal numbers in any
open interval, and we now describe that extension. Clearly we may
limit our attention to subintervals of $[0,1]$, since the set of
normal numbers to any base is invariant under translation by an
integer. In the original construction at the beginning of Section~2,
we began with $d_2=2^2$ and $\alpha_2=\frac34$; to be more general,
let $\alpha_2$ be any 2-adic fraction ${a\over d_2}$, where
$d_2=2^{n_2}$ for some positive integer $n_2$. Next fix any sequence
$n_3$, $n_4$, \dots of positive integers and modify the recursive
definition~(\ref{djdef}) of the $d_j$ to
\begin{equation}
d_j = j^{n_jd_{j-1}/(j-1)} \quad(j\ge3).  \label{newdjdef}
\end{equation}
If we now set
\begin{equation*}
\alpha_k = \alpha_2 \prod_{j=3}^k \big( 1-\frac1{d_j} \big)
\end{equation*}
(where of course the numbers $\alpha_k$ now depend on $\alpha_2$ and the
$n_j$), then the presence of the integers $n_j$ will not hinder the
proof of Lemma~\ref{intlem} that $d_k\alpha_k$ is always a $k$-adic
fraction (this is a consequence of the fact that $x-1$ always divides
$x^n-1$). Therefore the new limit $\alpha=\lim_{k\to\infty}\alpha_k$
can be shown to be a transcendental, absolutely abnormal number in
exactly the same way, the modifications only accelerating the
convergence of the infinite product and enhancing the ease with which
the various inequalities in the Lemmas are satisfied. (When we make this
modification, then the one case we must avoid is $\alpha_2=\frac12$
and $n_2=n_3=\dots=1$, for in this case it will happen that $d_j=j^1$
for every $j\ge2$. Then each product $\alpha_k$ is a telescoping
product with value $\frac1k$, and their limit $\alpha=0$, while
certainly absolutely abnormal, will be uninterestingly so.)

In particular, Lemma~\ref{squeezelem} applied with $k=2$ would show in
this context that ${a\over2^{n_2}}>\alpha>{a-2\over2^{n_2}}$; thus by choosing
$\alpha_2={a\over2^{n_2}}$ appropriately, we can ensure that the resulting
number $\alpha$ lies in any prescribed open subinterval of
$[0,1]$. Moreover, the various choices of the integers $n_3$, $n_4$,
\dots will give rise to distinct limits $\alpha$; one can show this by
considering the first index $j\ge3$ at which the choices of $n_j$
differ, say, and then applying Lemma~\ref{squeezelem} with $k=j$ to
each resulting $\alpha$. This generalization thus permits us to
construct uncountably many transcendental, absolutely abnormal numbers
in any prescribed open interval.

One interesting special case of this generalized construction arises
from the choices $\alpha_2=\frac12$ and $n_j=\phi(j-1)$ for all
$j\ge3$, where $\phi$ is the Euler totient function; these choices
give the simple recursive rule $d_j=j^{\phi(d_{j-1})}$ for $j\ge3$. In
this special case, the crucial property that $d_k\alpha_k$ is always
an integer is in fact a direct consequence of Euler's theorem that
$a^{\phi(n)}$ is always congruent to 1 modulo $n$ as long as $a$ and
$n$ have no common factors. In general, the smallest exponent $e_k$ we
can take in the recursive rule $d_k = k^{e_k}$ so that this crucial
property is satisfied is the multiplicative order of $k$ modulo $n$,
which might be smaller than $d_{k-1}/(k-1)$; however, our construction
given in Section~2 has the advantage of being more explicit, as it is
not necessary to wait and see the exact value of $d_{k-1}$ before
knowing how we will construct $d_k$.

We remark that Schmidt's construction \cite{schmidt} mentioned in the
introduction actually gives the following more powerful result: given
any set $S$ of integers exceeding 1 with the property that an integer $b$
is in $S$ if and only if every perfect power of $b$ is in $S$, Schmidt
constructed real numbers that are normal to every base $b\in S$ and
abnormal to every base $b\notin S$. (The problem considered herein is
the special case where $S=\emptyset$.) It would be interesting to see
if the construction in this paper could be modified to produce these
``selectively normal numbers'' as well.

We conclude with a few remarks about {\it absolutely simply abnormal
numbers\/}, numbers that are simply normal to no base whatsoever. As
we saw in the proof of our Theorem, the number $\alpha$ does in fact
meet this stronger criterion of abnormality. On the other hand, while
all rational numbers are absolutely abnormal, many of them are in fact
simply normal to various bases. For example, $1/3$ is simply normal to
the base 2, as its binary representation is $0.010101\dots$. In fact,
one can check that every fraction in reduced form whose denominator is
3, 5, 6, 9, 10, 11, 12, 13, 17, 18, 19, 20, \dots is simply normal to
the base 2 (presumably there are infinitely many such
odd denominators---can it be proven?); every fraction whose denominator is
7, 14, 19, 21, 31, \dots is simply normal to the base 3; and so
on. Somewhat generally, if $p$ is a prime such that one of the
divisors $b$ of $p-1$ is a primitive root for $p$, then every fraction
whose denominator is $p$ is simply normal to the base $b$ (although
this is not a necessary condition, as the normality of fractions with
denominator 17 to base 2 shows).

For a fraction with denominator $q$ to be simply normal to the base
$b$ (we can assume, by multiplying by $b$ a few times if necessary,
that $b$ and $q$ are relatively prime), it is necessary for $b$ to
divide the multiplicative order of $b$ modulo $q$, and hence $b$ must
certainly divide $\phi(q)$ by Euler's theorem. Therefore, we
immediately see that every fraction whose denominator is a power of 2
is absolutely simply abnormal. One can also verify by this criterion
that every fraction whose denominator in reduced form is 15 or 28, for
example, is absolutely simply abnormal. It seems to be a nontrivial
problem to classify, in general, which rational numbers are absolutely
simply abnormal. Since some fractions with denominator 63 are simply
normal to the base 2 while others are not, as one can check, absolute
simple abnormality probably depends in general on the numerator as
well as the denominator of the fraction.

\bigskip
{\smaller\smaller\narrower\noindent{\it Acknowledgements:\/} Glyn
Harman gave a survey talk on normal numbers at the Millennial
Conference on Number Theory at the University of Illinois in May 2000
(see~\cite{har}), at the end of which Andrew Granville asked the
question about a specific absolutely abnormal number that spurred this
paper. Carl Pomerance suggested the number $\sum_{n=1}^\infty
(n!)^{-n!}$, a Liouville number which again is almost surely
absolutely abnormal, but a proof of this seems hopeless. The author
would like to thank all three for their interest in this construction,
and also an anonymous referee for several helpful comments on a
preliminary version of this paper.\par}

\bibliography{abnormal}
\bibliographystyle{amsplain}
\end{document}